\documentclass[12pt]{article}
\usepackage{amsmath,amssymb}
\usepackage{geometry}
\usepackage{amsmath, amsfonts,amsthm,amssymb,enumerate,amscd}

\usepackage[active]{srcltx}
\usepackage[utf8]{inputenc}
\usepackage{algorithm}
\usepackage{stmaryrd}
\usepackage{array}
\usepackage{algorithmicx}
\usepackage{algpseudocode}
\usepackage{graphicx}
\usepackage{letltxmacro}

\usepackage{color}
\usepackage[table,xcdraw,svgnames,dvipsnames]{xcolor}
\definecolor{light-salmon}{RGB}{255,140,120}
\usepackage[colorlinks]{hyperref}
\hypersetup{
	allbordercolors=.,
	citecolor=DarkGreen,
	linkcolor=DarkRed,
	urlcolor=NavyBlue
}

\graphicspath{{./pics/}{./metapost/}}

\usepackage{enumitem}

\theoremstyle{plain}
\newtheorem{thm}{Theorem}
\newtheorem{lemma}{Lemma}

\theoremstyle{definition}
\newtheorem{rem}[thm]{Remark}

\newtheorem{prop}{Property}

\newcommand{\Per}{\operatorname{Per}}

\renewcommand{\Bbb}{\mathbb}
\newcommand{\bo}[1]{{\bf #1}}

\newcommand{\conv}{\operatorname{conv}}

\makeatletter
\DeclareFontFamily{U}{tipa}{}
\DeclareFontShape{U}{tipa}{m}{n}{<->tipa10}{}
\newcommand{\arc@char}{{\usefont{U}{tipa}{m}{n}\symbol{62}}}%

\newcommand{\arc}[1]{\mathpalette\arc@arc{#1}}

\newcommand{\arc@arc}[2]{%
	\sbox0{$\m@th#1#2$}%
	\vbox{
		\hbox{\resizebox{\wd0}{\height}{\arc@char}}
		\nointerlineskip
		\box0
	}%
}
\makeatother

\title{Mixed volumes and the Blaschke-Lebesgue theorem}

\author{Beniamin Bogosel}

\begin{document}
	
\maketitle

\begin{abstract}
	The mixed area of a Reuleaux polygon and its symmetric with respect to the origin is expressed in terms of the mixed area of two explicit polygons. This gives a geometric explanation of a classical proof due to Chakerian. Mixed areas and volumes are also used to reformulate the minimization of the volume under constant width constraint as isoperimetric problems. In the two dimensional case, the equivalent formulation is solved, providing another proof of the Blaschke-Lebesgue theorem. In the three dimensional case the proposed relaxed formulation involves the mean width, the area and inclusion constraints.
\end{abstract}

\section{Introduction}

Planar shapes of constant width have the remarkable property that the distance between any pair of enclosing parallel tangent lines is the same. It is not intuitive at first sight that there are many such shapes other than the disk. One classical example is the Reuleaux triangle: take an equilateral triangle of side length $1$ and consider the intersection of the three disks of radius $1$ centered at the vertices of the triangle. In the following $R$ denotes a Reuleaux triangle. Similar constructions can be made for regular $n$-gons having diameter $1$, obtaining the so-called Reuleaux regular $n$-gons. To fix the ideas, consider that all shapes have width $1$ in the following.

The constant width property does not allow such shapes to be too flat, therefore a valid question is to search for constant width shapes having minimal area. This question was answered by Blaschke \cite{Blaschke1915} and Lebesgue \cite{Lebesgue}:

{ 
\begin{thm}\label{thm:BL}
\bo{(Blaschke-Lebesgue)} The Reuleaux triangle is the planar shape of constant width having minimal area.	
\end{thm}
}

Various proofs of this result are known, some being summarized in \cite[Chapter 7]{yaglom-boltjanskii}. Nevertheless, there is a proof given in \cite{Chakerian} by Chakerian which seems too simple to be true. The main ideas are shown in Section \ref{sec:chakerian}, together with a geometric explanation of the surprising last step.

Let us recall that the Minkowski sum $A+B$ of two convex sets $A$ and $B$ is defined by
\[ A+B = \{a+b : a\in A , b \in B\}.\]
One should imagine translating a copy of $A$ with directions given by every element of $B$: $A+B = \bigcup_{b \in B}A+\{b\}$. Minkowski sums allow the following characterization for shapes of unit constant width:

\begin{prop} A convex body $K$ has unit constant width if and only if $K+(-K)=B$ where $B$ is the unit disk of radius $1$ centered at the origin.
	\label{prop:sum-cw}
\end{prop}
For a direct geometric proof see \cite[Exercise 7-15]{yaglom-boltjanskii}. Moreover, the following general result also holds and is attributed to P\`al \cite[p.127]{Eggleston}. 

\begin{prop}\label{prop:hexagon} Any plane convex shape of constant width $K$ can be inscribed in a regular hexagon of width $1$, i.e. a regular hexagon in which the distance between parallel opposite sides is equal to $1$.
\end{prop}

 A quick proof of this property is obtained observing that a rhombus with angles $\pi/3, 2\pi/3$ and width $1$ covers any orientation of $K$. Such a rhombus contains a unique unit width hexagon $H$. An elementary continuity argument shows that a rotation of $K$ inside the rhombus will be contained in the hexagon $H$.
 
\subsection{Structure of the paper} 
The objective of this paper is to illustrate how mixed volumes, introduced in the next section, provide solutions to Theorem \ref{thm:BL}. The original argument in Chakerian's proof \cite{Chakerian} is based on a series of estimates which turn out to be optimal. The key element of the proof is explained in Section \ref{sec:chakerian}. The mixed volume $A(K,-K)$ is constant when $K$ is contained between two particular polygons determined by a Reuleaux polygon. These facts are .

The \emph{skeleton} $\mathcal S(P)$ of a Reuleaux polygon $P$ is a polygon defined by the convex hull of its vertices. In addition to showing that the mixed area of $\mathcal S(P)$ and $-\mathcal S(P)$ coincides with the mixed area of $P$ and $-P$, the minimization of the area of $\mathcal S(P)$ is investigated in Theorem \ref{thm:skeleton-area}.

The convex hull $H = \conv\{R,-R\}$ of the Reuleaux triangle and its opposite plays an essential role in the proof given in \cite{Chakerian}.  Therefore, the study of the convex hull $S=\conv\{K,-K\}$ is motivated, for general shapes $K$ of constant width. Results of \cite{Gusev2010EulerCO}, \cite{Esterov2010TropicalVW}, \cite{schneider} allow to reformulate Theorem \ref{thm:BL} as an isoperimetric problem. A relaxation involving only convexity and inclusion constraints is proposed and solved in Section \ref{sec:chull}, giving a new proof for Theorem \ref{thm:BL}. 
 
The three dimensional case is discussed in Section \ref{sec:3D}. An obstacle to applying Chakerian's argument is the absence of a symmetric universal cover with sufficiently small volume in dimension three, the analogue of Property \ref{prop:hexagon}. Nevertheless, results in Section \ref{sec:chull} can be extended to the three dimensional case. The three dimensional Blascke-Lebesgue theorem is reformulated as an isoperimetric problem involving the surface area and the mean width. Hopefully, this will provide a new strategy for approaching the three dimensional analogue of Theorem \ref{thm:BL}. 

The illustrations in this paper are produced using Metapost and Matlab.

\section{Geometric ideas about Chakerian's proof}
\label{sec:chakerian}

\subsection{Mixed areas}

The mixed area of two planar shapes of constant width is defined using the following relation \cite[Chapter 5]{schneider}
\begin{equation}\label{eq:mixed}
 A(K_1+K_2) = A(K_1)+2A(K_1,K_2)+A(K_2).
\end{equation}
The mixed area is symmetric and monotone in each one of the arguments: $K_1\subset K_2 \Longrightarrow A(K_1,K)\leq A(K_2,K)$. Moreover, since $K+K=2K$ it follows from \eqref{eq:mixed} that $A(K,K) = A(K)$, the area of the convex set $K$.

Taking $K_1=K, K_2=-K$ in \eqref{eq:mixed} and using Property \ref{prop:sum-cw} we find that
\begin{equation}\label{eq:cw}
 \pi = A(B)=  A(K+(-K)) = 2A(K)+2A(K,-K).
\end{equation}
Therefore, denoting with $\mathcal{CW}$ the class of planar shapes with unit constant width, the following problems are equivalent:
\begin{equation}
\min_{ K \in \mathcal{CW}} A(K) \ \ \Longleftrightarrow \max_{K \in \mathcal{CW}} A(K,-K).
\end{equation}
This motivates to search for upper bounds of $A(K,-K)$ using the monotonicity of the mixed area with respect to inclusion.
Property \ref{prop:hexagon} states that there exists a unit width regular hexagon $H$ such that $K\subset H$. Assuming $H$ is centered at the origin, we also have $-K \subset H$ so by monotonicity of the mixed volume we have
\[ A(K,-K) \leq A(K,H)\leq A(H,H) = A(H) = \sqrt{3}/2.\]
Thus, using \eqref{eq:cw} we find that any constant width shape $K$ verifies
\[ A(K) = \frac{\pi}{2}-A(K,-K) \geq \frac{\pi}{2}-A(H) =\frac{\pi-\sqrt{3}}{2}=A(R).\]
The lower for the area $A(K)$ found above coincides with the area of the Reuleaux triangle. Thus, the Blaschke-Lebesgue theorem is proved. 

\begin{figure}
	\centering
	\begin{tabular}{ccc}
	\includegraphics[height=0.32\textwidth]{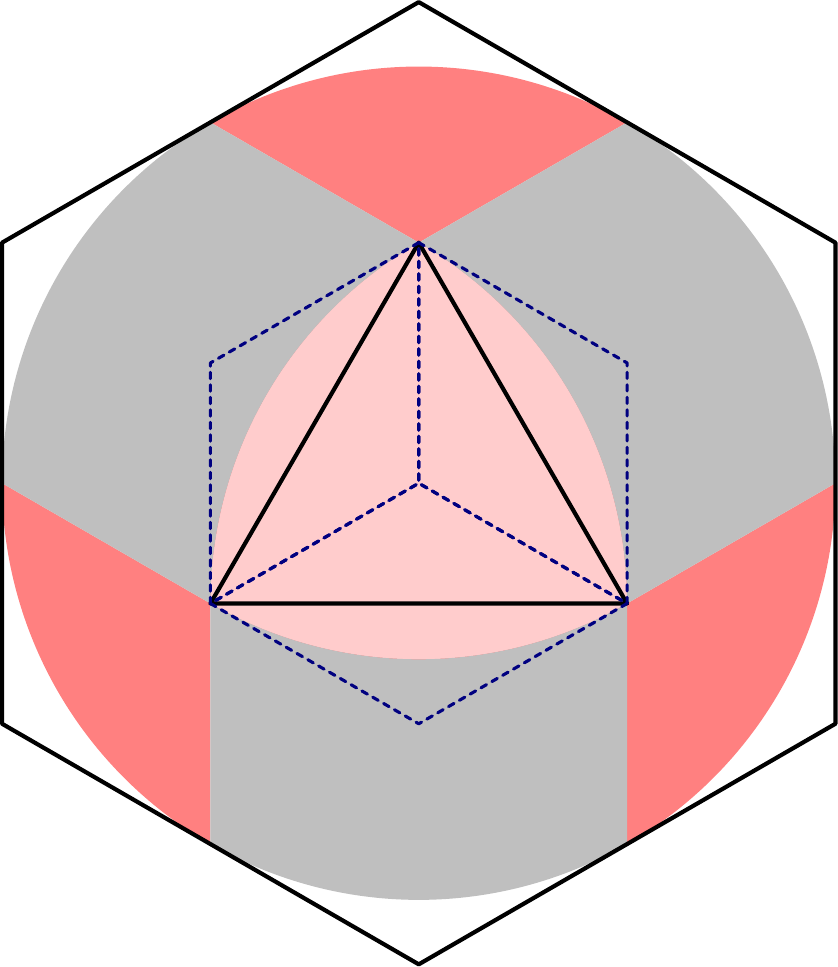}&
	\includegraphics[height=0.32\textwidth]{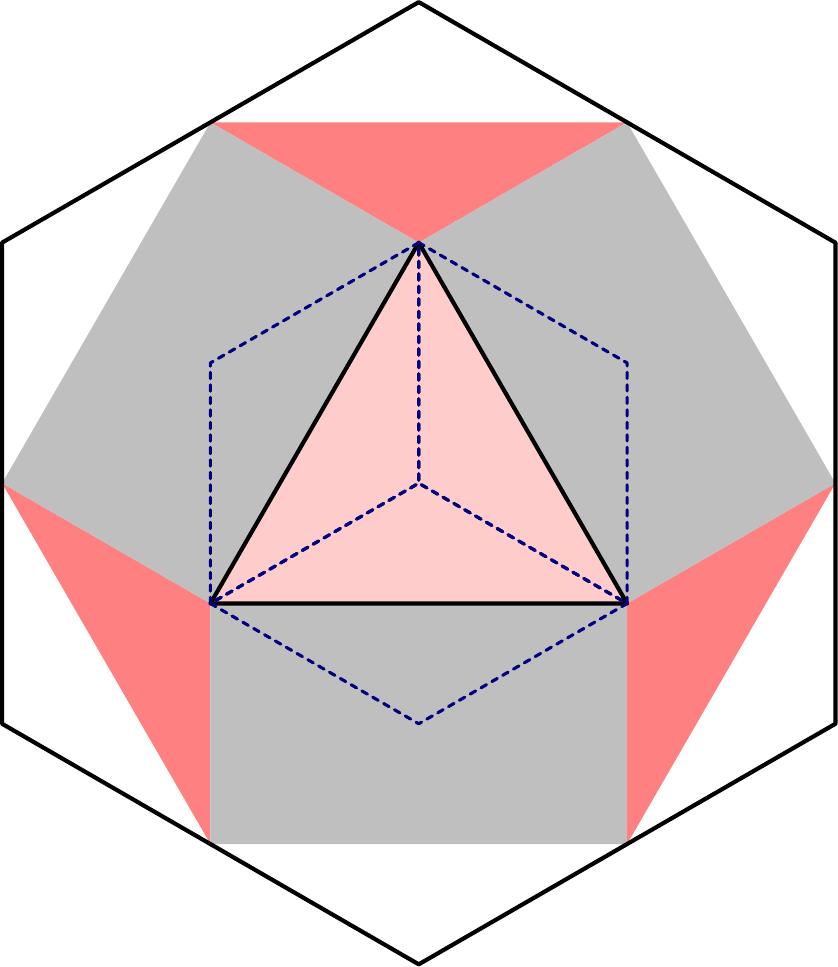}&
	\includegraphics[height=0.32\textwidth]{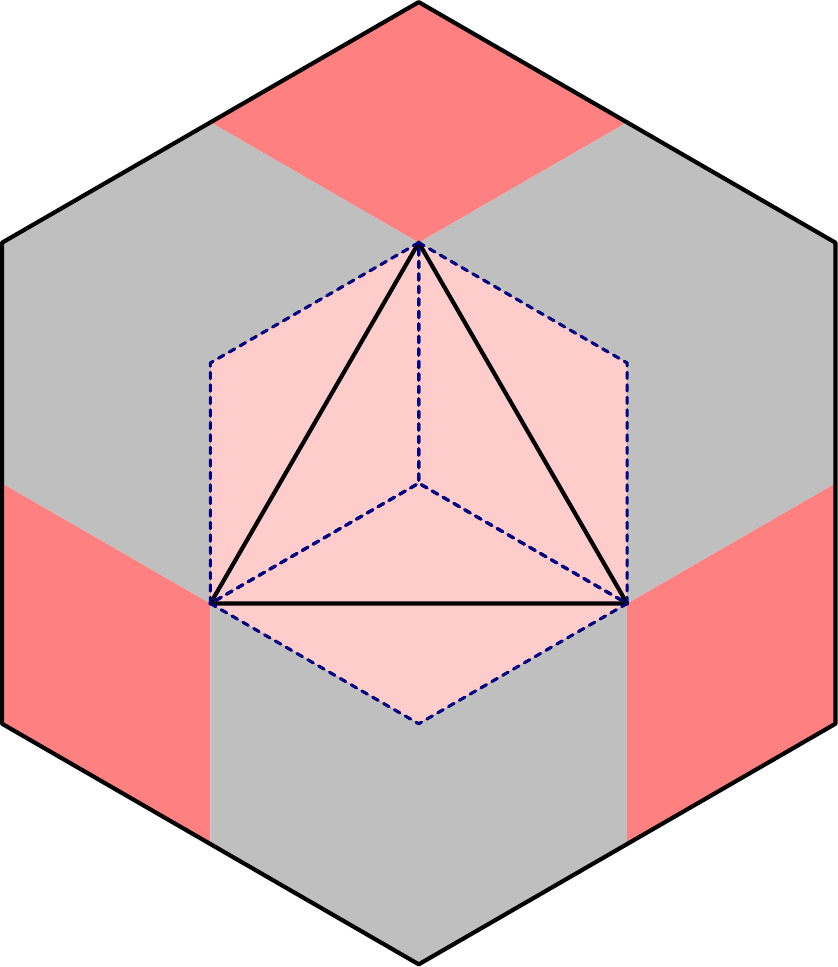}\\
	$R+(-R)$ & $T+(-T)$ & $H+(-H)$
	\end{tabular}
\caption{Computation of the Minkowski sum $K+(-K)$ for the Reuleaux triangle, its inscribed equilateral triangle and its circumscribed regular hexagon. The central area represents $K$, the red outer slices are thirds of $-K$ while the grey areas are have equal $2A(K,-K)$. It can be observed that the grey areas are the same for the three figures.}
\label{fig:graphical-mixed-volume}
\end{figure}

It is surprising that replacing $K$ and $-K$ by the covering hexagon of width $1$ provides a tight enough inequality to prove the result! In particular, it follows that $A(-R,R) = A(-H,H) = A(H)$, showing that mixed volumes are not necessarily strictly increasing in their arguments.

Indeed, in Figure \ref{fig:graphical-mixed-volume} three configurations are shown, computing the Minkowski sum $K+(-K)$ for the Reuleaux triangle, its inscribed equilateral triangle and the circumscribed regular hexagon. Splitting $-K$ into three slices of central angle $2\pi/3$, three types of areas are apparent in the figures:

(a) The central one in pink corresponds to the original shape $K$.

(b) The outer slices shown in red correspond to the opposite shape $-K$. Their total area equals $A(K)$.

\begin{figure}
	\includegraphics[width=\textwidth]{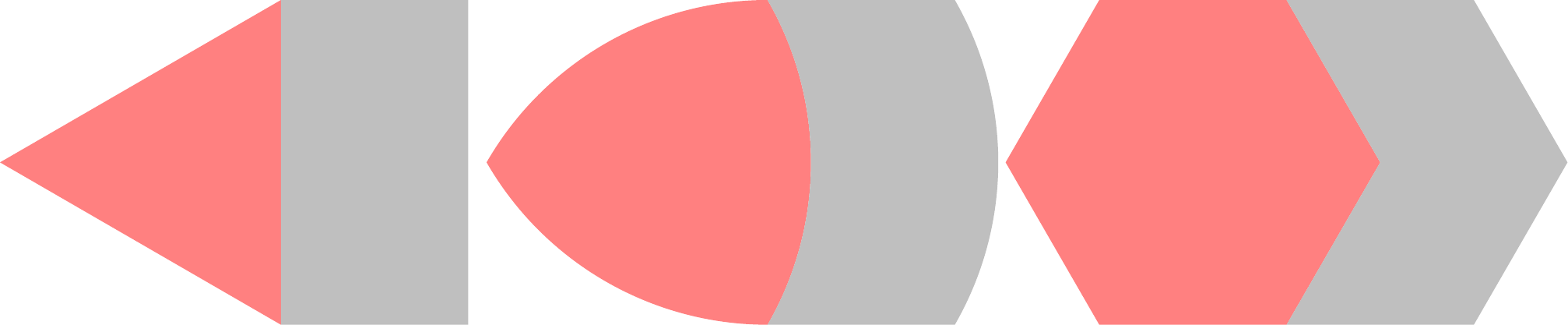}
	\caption{Drawing of the three regions corresponding to $2A(K,-K)$ in Figure \ref{fig:graphical-mixed-volume} for the equilateral triangle, the Reuleaux triangle and the regular hexagon. The three gray regions have equal heights and horizontal slices with equal constant lengths. Thus they have the same area.}
	\label{fig:parallel-slices}
\end{figure}

(c) The grey areas, according to \eqref{eq:mixed}, correspond to the mixed volume $2A(K,-K)$. It can be observed that the grey areas are the same for the three cases since each one of the three regions has height $1$ and horizontal slices equal to $\sqrt{3}/{3}$, the circumradius of the unit equilateral triangle. See Figure \ref{fig:parallel-slices} for an illustration.

Thus, $A(T,-T)=A(R,-R)=A(H,-H)=A(H)$ and more generally if $T\subset K \subset H$ then $A(K,-K) = A(H)$. This surprising result will be extended below to all Reuleaux polygons. 

\subsection{A general property regarding Reuleaux polygons}

\begin{figure}
	\centering 
	\includegraphics[height=0.3\textwidth]{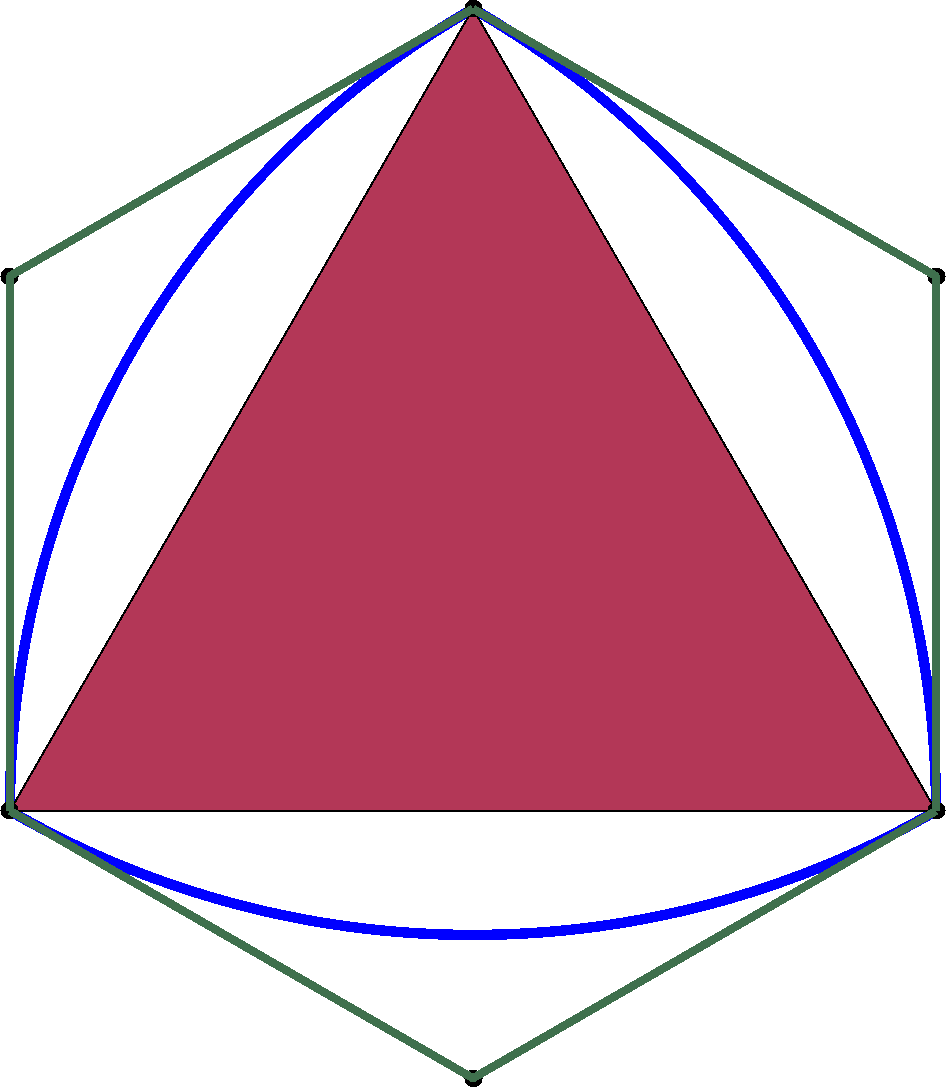}\qquad
	\includegraphics[height=0.3\textwidth]{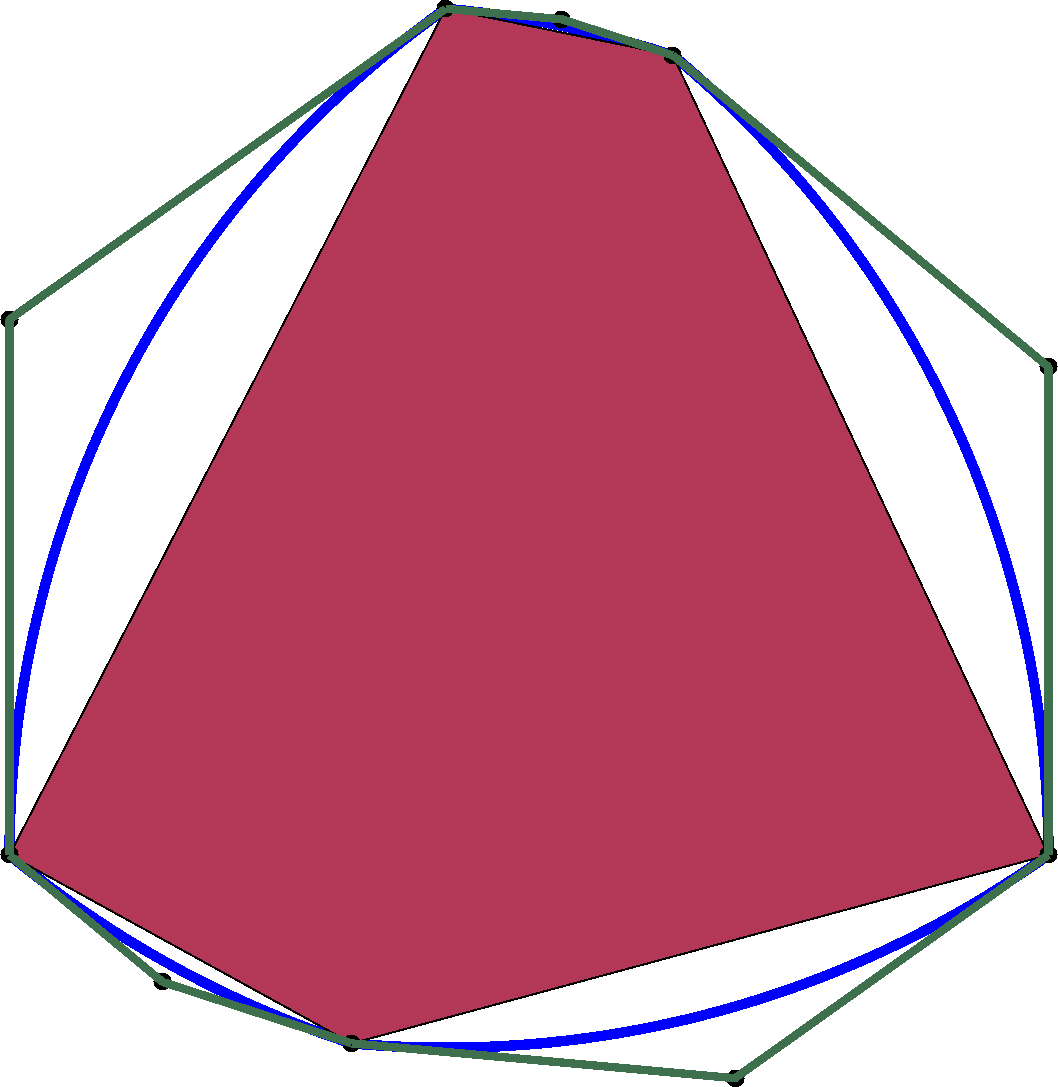}\qquad
	\includegraphics[height=0.3\textwidth]{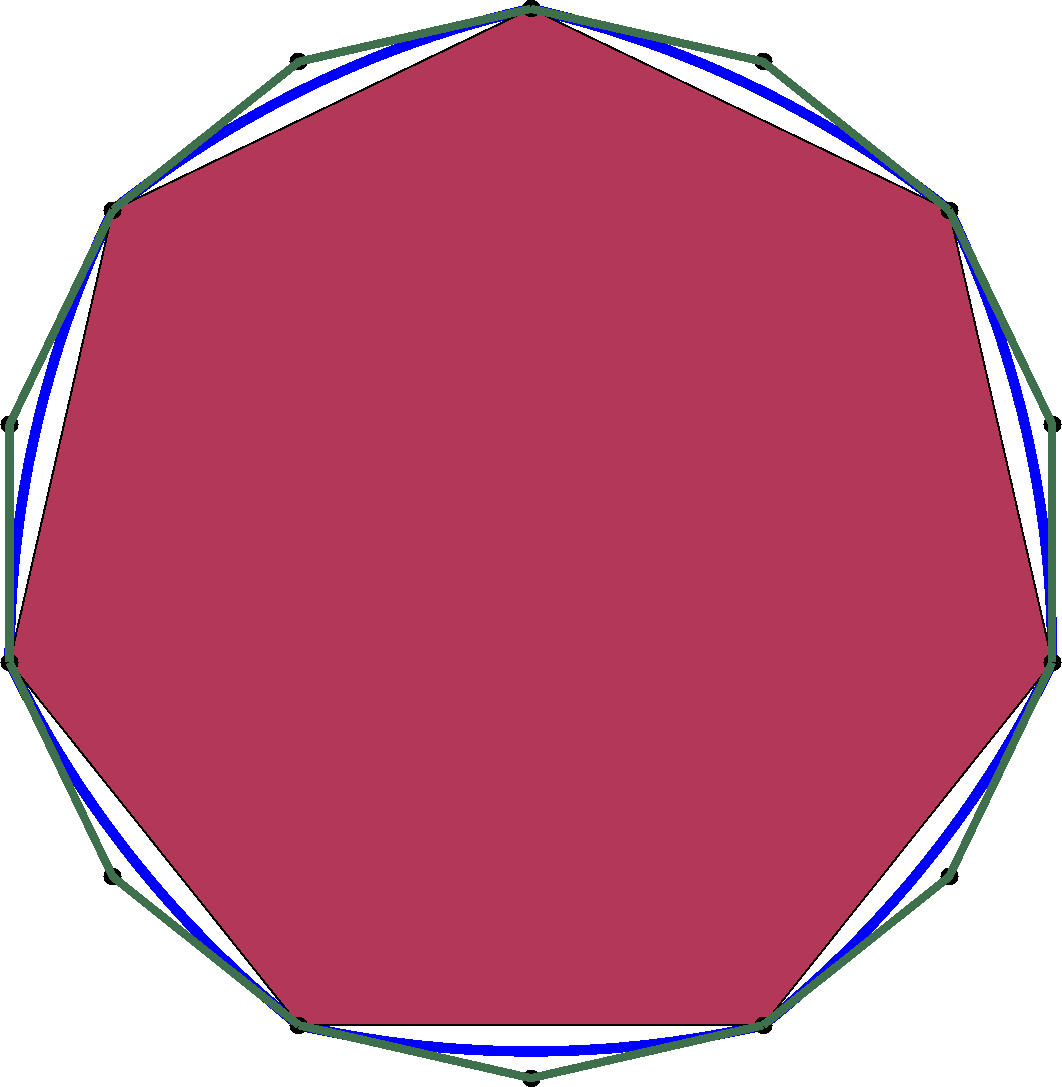}
	\caption{Examples of Reuleaux polygons $P$, their associated skeleton polygons $\mathcal S(P)$ and their tangent circumscribed polygons $\mathcal C(P)$. The mixed area $A(K,-K)$ is constant for $\mathcal S(P) \subset K \subset \mathcal C(P)$.} 
	\label{fig:skeleton}
\end{figure}

The property underlined above and described in Figure \ref{fig:graphical-mixed-volume} also extends to general Reuleaux polygons. Given a general Reuleaux polygon $P$ with $n=2k+1 \geq 3$ sides, consider the following elements, illustrated in Figure \ref{fig:skeleton}:

\begin{itemize}[noitemsep]
	\item $\mathcal S(P)$ is the \emph{skeleton polygon} of $P$, obtained by removing all circle caps from $P$. Equivalently, $\mathcal S(P)$ is the convex hull of the vertices of $P$.
	\item $\mathcal C(P)$ is the \emph{tangent circumscribed polygon} of $P$, obtained by drawing extremal tangent lines to all arcs in the boundary of $P$. 
\end{itemize}

\begin{thm}\label{thm:mixed-constant}
   If $P$ is a Reuleaux polygon then 
   \[ A(K,-K) = A(P,-P)=A(\mathcal S(P),-\mathcal S(P))=A(\mathcal C(P),-\mathcal C(P)),\] 
   for all $K$ which verify $\mathcal S(P) \subset K \subset \mathcal C(P)$. 
\end{thm}

\emph{Proof:} Consider $x$ a vertex of $P$. Denote by $\bo n(x)$ the normal vectors at $x$. Obviously, if $n \in \bo n(x)$ then $-n$ is a normal vector for the arc opposite to $x$ in $\partial P$. Denote by $\alpha(x)$ the angle with vertex at $x$ determined by the opposite arc in $\partial P$. Any convex shape verifying $\mathcal S(P) \subset K \subset \mathcal C(P)$ has the following properties:
\begin{itemize}[noitemsep]
	\item $x$ is also a vertex of $K$ and $\bo n(x)$ is contained in the set of normal vectors to $\partial K$ at $x$.
	\item normal vectors to $\partial K \cap \alpha(x)$ are given by $-\bo n(x)$.
\end{itemize}
Denoting with $o$ the origin, inside $P$, it follows that in the Minkowski sum $K+(-K)$, has the following boundary structure, for each vertex $x$ of $P$:
\begin{itemize}[noitemsep]
	\item Reflect $\partial K \cap \alpha(x)$ about the origin and translate it with the vector $\overrightarrow{ox}$
	\item Translate $\partial K \cap \alpha(k)$ with vector $\overrightarrow{xo}$.
\end{itemize}
In particular, $K+(-K)$ can be partitioned into two copies of $K$ and a series of curved parallelograms with sides $(\partial K \cap \alpha(x)) \times ox$ (configuration similar to Figure \ref{fig:parallel-slices}). The area of each such parallelogram does not depend on $K$, once the origin is fixed. For each vertex $x$ of $P$, the width of the parallelograms is equal to the length of the projection of $\partial K \cap \alpha(x)$ on the line orthogonal to $ox$ and the horizontal sections all have the same length $|ox|$.  Since $2A(K,-K) = A(K+(-K))-2A(K)$ it follows that $2A(K,-K)$ is the total area of the curved parallelograms, which does not depend on $K$. \hfill $\square$

The proof of Theorem \ref{thm:mixed-constant} can also be achieved using the valuation property of mixed volumes \cite[Chapter 5]{schneider}. If $K_1,K_2$ are convex such that $K_1\cup K_2$ is convex then
\[ A(K_1\cup K_2,K)+A(K_1\cap K_2,K) = A(K_1,K)+A(K_2,K).\]
The valuation property is not straightforward to use, since even when splitting a convex body along a segment, the mixed volume involving that segment cannot be neglected. 

\noindent \bo{Example.} Imagine partitioning the unit disk $B$ into the two half disks $B_1,B_2$. Recall that the two dimensional steiner formula gives $A(K,B) = \frac{1}{2}P(B)$. Then 
\[ A(B,B_1)+A(B,B_2) = A(B,B)+A(S,B),\]
where $S$ is the segment $B_1\cap B_2$. Note that $A(S,B) = 1>0$ even though $S$ has zero area. 

Results in the following Lemma can quickly be obtained using basic properties of Minkowski addition and mixed volumes.

\begin{lemma}\label{lem:segment}
	(a) Let $K$ be a convex shape and $S$ a segment of length $l$. Then $A(K,S) = \frac{1}{2} l w$ where $w$ is the width of $K$ in the direction orthogonal to $S$.
	
	(b) Let $L$ be a shape whose boundary contains the segment $S$. Suppose that $K$ is such that there exists $x \in \partial K$ and any normal vector to $\partial L\setminus S$ is a normal vector of $K$ at $x$. Then 
	\[ A(L,K) = A(S,K).\]
\end{lemma}

\emph{Proof:} (a) It is obvious that $K+S$ widens $K$ with the length of $S$ in the direction of $S$. A curved band of height $w$ having slices with constant length $l$ is added to $K$. The result follows at once from \eqref{eq:mixed} and $A(S)=0$.

(b) The property concerning normal vectors of $\partial L \setminus S$ and the normals of $K$ at $x$ imply that $K+L$ is the union of $K+S$ and a translate of $L$. The result follows from \eqref{eq:mixed}. \hfill $\square$

\begin{rem}
	It is possible to obtain Theorem \ref{thm:mixed-constant} by repeatedly applying Lemma \ref{lem:segment} for the circular segments contained in a Reuleaux polygon $P$. Denote by $K$ the skeleton of $P$ and let $L_1,...,L_n$ be the circular segments near each arc of the boundary of $P$. Denote by $S_i$ the segment in the boundary of $L_i$. Then $-K$ and each one of the $L_i, S_i$, $i=1,...,n$ verify the hypothesis of Lemma \ref{lem:segment} part (b). Then we have
	\[ A(K\cup L_i,-K) = A(K,-K)+A(L_i,-K)-A(S_i,-K) = A(K,-K),\]
	since $A(L_i,-K)=A(S_i,-K)$. 
	The argument can be repeated for each side of $K$.
\end{rem}

Skeletons of Reuleaux polygons behave in a similar way with respect to area maximization and minimization. We have the following consequence of the Blaschke-Lebesgue theorem.

\begin{thm}\label{thm:skeleton-area}
	If $\mathcal S(P)$ is the skeleton of a Reuleaux polygon $P$ then the area of $\mathcal S(P)$ is minimal when $P$ is a Reuleaux triangle and $\mathcal S(P)$ is an equilateral triangle. The area of $\mathcal S(P)$ is bounded above by $\pi/4$, the area of a disk of radius $1/2$ and the extremal value is not attained.
\end{thm}

For the upper bound, there is nothing to prove. A skeleton $\mathcal S(P)$ is strictly contained in $P$. A Reuleaux polygon has area strictly bounded above by $\pi/4$, consequence of the isodiametric inequality or the Alexandrov-Fenchel inequality for mixed volumes $A(K_1,K_2)^2 \geq A(K_1)A(K_2)$ (with equality only when $K_1,K_2$ are homothetic). 

For the lower bound two proofs are given. The first one is a consequence of the Blaschke-Lebesgue theorem. The second one is a consequence of the inequality
\begin{equation}\label{eq:diff-body-2D}
 A(K+(-K)) \leq 6 A(K) \Longleftrightarrow A(K,-K) \leq 2A(K),
 \end{equation}
verified for all two dimensional convex shapes. See \cite{rogers-shepard} for details and other references. Equality holds in \eqref{eq:diff-body-2D} if and only if $K$ is a simplex.

\emph{Proof:} {\bf (a) Using the Blaschke-Lebesgue theorem.} Denoting by $\theta_1,...,\theta_n$ the angles which determine the edges of $P$, we obviously have $\theta_1+...+\theta_n =\pi$ and $\theta_i \in [0,\pi/3]$ (see \cite[Chapter 7]{yaglom-boltjanskii} for the proofs). Denote by $f(\theta) = (\theta-\sin \theta)/2$, the area of a circular segment with angle $\theta$. Then 
\[ A(\mathcal S(P)) = A(P)-\sum_{i=1}^n f(\theta_i).\]
Differentiating $f$ twice, it is immediate to deduce that $f$ is convex. Therefore, the maximal value of $\sum_{i=1}^n f(\theta_i)$ is $3f(\pi/3)$ under the given constraints for the angles. Since the area is also minimized by the Reuleaux triangle $R$, it follows that
\[ A(\mathcal S(P)) \geq A(R)-3f(\pi/3) = A(\mathcal S(R)).\]

{\bf (b) Using \eqref{eq:diff-body-2D}.} Keep the notation above for the area of a circular segment. The proof of Theorem \ref{thm:mixed-constant} implies that if $P$ is a Reuleaux polygon with angles $\theta_1,...,\theta_n$ then 
\[ \pi = A(P+(-P)) = A(\mathcal S(P)+(-\mathcal S(P)))+2\sum_{i=1}^n f(\theta_i).\]
This implies, in particular, that $A(\mathcal S(P)+(-\mathcal S(P))) = \sum_{i=1}^n \sin \theta_i$.

Using the fact that $\sum_{i=1}^n f(\theta_i) \leq 3f(\pi/3)$ and \eqref{eq:diff-body-2D} we find that
\[ \pi \leq 6A(\mathcal S(P))+6f(\pi/3) = 6A(\mathcal S(P))+\pi-3\sqrt{3}/2.\]
Thus $A(\mathcal S(P))\geq \sqrt{3}/4$, the area of an equilateral triangle with unit width. This coincides with the area of the skeleton of the Reuleaux triangle, implying the result. Equality holds in the estimates above only if $P$ is the Reuleaux triangle $R$. \hfill $\square$

Denote by $\mathcal S$ the family of skeletons of Reuleaux polygons of unit width. The Blaschke Lebesgue theorem and Theorems \ref{thm:mixed-constant}, \ref{thm:skeleton-area} show that 
\begin{equation}\label{eq:duality}
 \max_{Q \in \mathcal S}A(Q,-Q) \leq A(T,-T) = 2A(T) \leq \min_{\Omega \in \mathcal S}2A(\Omega).
\end{equation}
Therefore, a duality-like relation holds for the maximization of the mixed area $A(K,-K)$ and the minimization of the area in the class of skeletons of Reuleaux polygons.

Moreover, if \eqref{eq:duality} is proved independent from Theorem \ref{thm:BL}, then it implies the Blaschke-Lebesgue theorem, since for any Reuleaux polygon, according to Theorem \ref{thm:mixed-constant}, we have
\[ A(P,-P) = A(\mathcal S(P),-\mathcal S(P)) \leq 2A(T) = A(R,-R).\] 

\noindent\bo{Universal covers.} The arguments above also imply the following result.

\begin{thm}\label{eq:thm-sym-univ-cover}
	The central symmetric shape having minimal area which covers a copy of any set of diameter at most $1$ is the regular hexagon of unit width.
\end{thm}

The general universal cover problem in the plane is still unsolved and recent developments can be found in \cite{UnivCoverBaez}, \cite{gibbs2018upper}. If the symmetry condition is removed, various slices of the regular hexagon around particular corners can be removed. 

In dimension two, the best symmetric universal cover $H$ is the convex hull of $R$ and $-R$. A surprising formula for computing the area of the convex hull of $K$ and $-K$ leads to a new proof of the two dimensional Blaschke-Lebesgue theorem, presented in the following section.

\section{Reverse isoperimetric inequalities}
\label{sec:chull}

The following relation, computing the area of the convex hull of two convex shapes, was discovered in \cite{Gusev2010EulerCO} and \cite{Esterov2010TropicalVW} in the context of tropical geometry and extended by Schneider in \cite{Schneider_mix_vol}. In dimension two, if $S_1,S_2$ are two convex bodies and their convex hull is $S = \conv\{S_1,S_2\}$, then we have
\[ A(S)-A(S,S_1)-A(S,S_2)+A(S_1,S_2) = 0.\]
In particular, if $S = \conv \{K,-K\}$ we obtain
\begin{equation}\label{eq:reformulation} A(K,-K) = A(S,K)+A(S,-K)-A(S) = A(S,B)-A(S) = \frac{1}{2}\Per(S)-A(S),
\end{equation}
where $\Per(S)$ is the perimeter of the convex set $S$.
Therefore, we obtain the following result.

\begin{thm}\label{thm:reformulation} Finding the shape of constant width with minimal area is equivalent to solving the problem
\begin{equation}
\max_{S = \conv\{K,-K\}} \frac{1}{2}\Per(S)-A(S).
\label{eq:max_equiv_2D}
\end{equation}
where $S$ is the convex hull of a shape of constant width $K$ and its symmetric with respect to the origin $-K$.
\end{thm} 

The equivalent problem is not necessarily simpler than the original one. Nevertheless, solutions to problems of this kind are completely characterized in \cite{bianchini_henrot} in the case where $S$ is contained in an annulus.

Denoting $B(r)$ the disk/ball of radius $r$ centered at the origin, it is not difficult to see that $S=\conv\{K,-K\}$ is contained (up to a translation) in the annulus determined by the disks $B(\sqrt{3}/{3})$ and $B(1/2)$. First, observe that according to \eqref{eq:reformulation} the objective function in \eqref{eq:max_equiv_2D} is invariant to translations of $K$ and $-K$. Without loss of generality, suppose their circumcenters coincide. Then we have the following:

\begin{itemize}[noitemsep]
	\item The circumradius of a two dimensional shape of unit constant width is bounded above by $\sqrt{3}/3$. See \cite[Exercise 7-14]{yaglom-boltjanskii}, for example.
	\item Given a direction $\theta \in [0,\pi]$, the distance from the origin to the tangent plane to $K$ orthogonal to $\theta$ is the support function $h_K$ of $K$. It is not difficult to observe \cite[p. 140]{schneider} that if $K$ has constant width one then $h_K(\theta)+h_K(\theta+\pi)=1$ for all $\theta \in [0,2\pi]$. Moreover, the symmetric body verifies $h_{-K}(\theta) = h_K(\theta+\pi)=1-h_K(\theta)$. The convex hull of $K,-K$ verifies $h_S(\theta) = \max \{ h_K(\theta), 1-h_K(\theta)\}$. Thus, in any case, $h_S(\theta) \geq 1/2$, implying that $S$ contains the disk of radius $1/2$.
\end{itemize}

Denoting $\mathcal A$ the family of convex shapes $K$ verifying $B(1/2)\subset K \subset B(\sqrt{3}/3)$, consider the following relaxation of \eqref{eq:max_equiv_2D}:

\begin{equation}
\max_{S \in \mathcal  A} \frac{1}{2}\Per(S)-A(S).
\label{eq:max_equiv_2D_relax}
\end{equation}

It is equivalent to minimize $2A(S)-\Per(S)$ for $S\in \mathcal A$. A quick investigation shows that \cite[Theorem 3.1]{bianchini_henrot} applies for $a=\frac{1}{2}$, $b=\frac{\sqrt{3}}{2}$, $\lambda=2$ showing that \eqref{eq:max_equiv_2D_relax} is solved by the regular hexagon inscribed in $B(\sqrt{3}/3)$ and circumscribed to $B(1/2)$. This hexagon has edge length $\sqrt{3}/3$, its perimeter is $\Per(H) = 2\sqrt{3}$ and its area is $A(H)=\sqrt{3}/2$. The maximal value attained in \eqref{eq:max_equiv_2D_relax} is $\sqrt{3}/2$.

The regular hexagon inscribed in $B(\sqrt{3}/3)$ is also admissible in \eqref{eq:max_equiv_2D}, since $H$ is the convex hull of $R\cup (-R)$, the Reuleaux triangle and its symmetric with respect to the origin. Therefore, the regular hexagon is also a solution for \eqref{eq:max_equiv_2D}, providing another proof of the Blaschke-Lebesgue theorem is obtained. In addition, a new strategy for attacking the three dimensional case is provided.

\begin{rem}
	It is not possible to further relax problem \eqref{eq:max_equiv_2D_relax}, using only inclusion constraints, without changing the solution. Considering only sets $S$ with $B(1/2)\subset S$ allows the unit square as an admissible shape, with objective value equal to $1$, strictly greater than $\sqrt{3}/2$. Considering only the constraint $S\subset B(\sqrt{3}/3)$ and a sequence of rectangles converging to a diameter leads to an objective value of $2\sqrt{3}/3$, strictly larger than $\sqrt{3}/2$.
\end{rem}

\section{Discussion around the three dimensional case}
\label{sec:3D}

The body of constant width and minimal volume is not yet completely characterized in dimension three. Let us first underline that such a body $K^*$ of constant width and minimal volume exists. This follows from the classical Blaschke selection theorem. Historical facts about the problem are gathered in \cite{kawohl-webe}. In particular, the area and volume of a constant width body $K$ in dimension three are linked through an explicit formula due to Blaschke: $V(K)=\frac{1}{2}\Per( K)-\frac{\pi}{3}$, where $V(K)$, $\Per(K)$ denote the volume and surface area of the three dimensional convex body $K$. Therefore, minimizing the volume is equivalent to minimizing the area.

The interest in the three dimensional analogue of Theorem \ref{thm:BL} was rekindled recently by a series of papers. In \cite{Hynd_Doubly_Monotone} a flow is introduced in the class of three dimensional bodies of constant width, which simultaneously increases the volume and decreases the circumradius. Reversing such a flow would have a great impact on understanding the three dimensional case. The three dimensional analogue of Reuleaux polygons were introduced in \cite{montejano}. Their structure was investigated in detail in \cite{meissner_hynd} and their surface perimeters were computed in \cite{hynd-vol-per} and \cite{bogosel_Meissner}. In particular, in \cite{bogosel_Meissner} the three dimensional case is reduced to a series of explicit finite dimensional problems.

It is conjectured that the three dimensional analogue of the Reuleaux triangle is the constant width shape minimizing the volume. The construction of these tetrahedra, called Meissner bodies, is shown in \cite[Chapter 7]{yaglom-boltjanskii}. The intersection of four balls centered at the vertices of a regular tetrahedron gives the Reuleaux tetrahedron, which is not a shape of constant width. For all three pairs of opposite edges, one edge should be smoothed to recover a constant width shape. In \cite[Section 8.3]{bodies_of_constant_width} more details can be found.

The best analytical lower bound for the volume of $K^*$ is also due to Chakerian \cite{Chakerian}. This follows from the classical inequality $3V(K) \geq r \Per(K)$ where $r$ is the inradius of $K$. A short proof using the support function of $K$ is contained in \cite[eq. (10)]{Chakerian}. It is known that $r \geq 1-\sqrt{3/8}$ for three dimensional bodies of constant width \cite[Section 14.3]{bodies_of_constant_width}. The desired bound follows, using the Blaschke formula:
\begin{equation}\label{eq:ineq-vol-inradius}
 V(B) \geq \frac{2r\pi}{3(3-2r)}\geq \frac{\pi}{3}(3\sqrt{6}-7)\approx 0.365.
\end{equation}
If $V_M = \pi\left(\frac{2}{3}-\frac{\sqrt{3}}{4} \arccos \frac{1}{3}\right)$ denotes the volume of a Meissner tetrahedron, solving the inequation 
\[  \frac{2r\pi}{3(3-2r)}\geq V_M\]
shows that constant width bodies with inradius
\[ r \geq r_M= \frac{3V_M}{\frac{2\pi}{3}+2V_M} \approx 0.429\]
do not minimize the volume. Therefore, we obtain the following.
\begin{thm}
The inradius of a unit three dimensional constant width body with minimal volume verifies
\[ r \in \left[1-\sqrt{\frac{3}{8}},r_M\right].\]
The circumradius $R = 1-r$ verifies $R \in [1-r_M, \sqrt{3/8}]$.
\end{thm}

 While mixed volumes are also defined in dimension three, the missing ingredient for copying Chakerian's proof in higher dimensions is knowing a good enough universal cover in dimension three. In \cite{kuperberg} it is shown that the rhombic dodecahedron $D$ having unit width is a universal cover. The length of the edge of $D$ is $a = \sqrt{3/8}$. Its volume is $\frac{16\sqrt{3}}{9}a^3$ and its area is $8\sqrt{2}a^2$.

Three dimensional mixed volume generates a trilinear form $V$. A three dimensional constant width body $K$ also verifies $K+(-K)=B$, therefore
\begin{equation}\label{eq:KplusK3D} \frac{4\pi}{3}= V(B) = V(B,K+(-K),K+(-K)) =2 V(B,K,K)+2V(B,K,-K).\end{equation}
and
\begin{equation}\label{eq:KplusK3Dvol}
\frac{4\pi}{3} = V(B) = 2V(K)+6V(K,K,-K).
\end{equation}
Like in the two dimensional case, denoting by $\mathcal{CW}$ the class of bodies of unit constant width, the following problems are equivalent
\[ \min_{K \in \mathcal{CW}} V(K) \Longleftrightarrow \max_{K \in \mathcal{CW}} A(K,K,-K) \Longleftrightarrow \max_{K \in \mathcal{CW}} A(B,K,-K).\]

Using $D$ as an universal cover gives
\[ V(K,K,-K)\leq V(D,D,D) = V(D) = \sqrt{2}/2.\]
However, combined with \eqref{eq:KplusK3Dvol} this gives a meaningless inequality:
\[ 2V(K) \geq \frac{4\pi}{3}-6\frac{\sqrt{2}}{2}=-0.053... \] 

Elementary properties of the mixed volume show that for $\delta>0$ we have
\[ V(K+\delta B) = V(K)+3\delta V(B,K,K)+3\delta^2 V(B,B,K)+\delta^3 V(B).\]
On the other hand, the classical  Steiner formula \cite{steiner_formula} gives
\[ V(K+\delta B) = V(K)+\Per(K)\delta +\frac{1}{2}k_1(\partial B)\delta^2 + \delta^3 V(B),\]
where for a smooth shape $B$, $k_1(B)$ denotes the average mean curvature. Direct identification shows that
\[ V(B,K,K) = \frac{1}{3}\Per(K),\ \ V(B,B,K) = \frac{1}{6}k_1(B).\]
In dimension three the average mean curvature is related to the mean width $\mathcal M(K)$ defined as follows. Given a direction $\theta \in \Bbb{S}^2$ compute the \emph{breadth} $b(K,\theta)$ of $K$ in the direction $\theta$ as the distance between the two parallel tangent planes having normal vector $\theta$ which enclose $K$. The mean width is just the average of the breadth for all possible directions
\[ \mathcal M(K) = \frac{1}{\Per(B)} \int_{\Bbb{S}^2} b(\theta,K)d\sigma.\]
In \cite[Chapter 9]{TopicsConvex} it is shown that 
\begin{equation}\label{eq:meanw-mixed-vol}
V(B,B,K) = \frac{2\pi}{3} \mathcal M(K).
\end{equation}

In view of \eqref{eq:KplusK3D} it follows that the Blaschke-Lebesgue problem in dimension three is equivalent to maximizing $V(B,K,-K)$ when $K$ has constant width. 
Like in the previous section, \cite{Schneider_mix_vol} implies that if $S=\conv\{K,-K\}$ and $K$ has constant width then 
\begin{equation}\label{eq:schneider}
V(B,K,-K) =  V(B,S,K)+V(B,S,-K)-V(B,S,S) = V(B,B,S)-V(B,S,S).
\end{equation}
We obtain the following.
\begin{thm}\label{thm:reformulation3D} Finding the shape of constant width with minimal area is equivalent to solving the problem
	\begin{equation}
	\max_{S = \conv\{K,-K\}} 2\pi \mathcal M(S)-\Per(S).
	\label{eq:max_equiv_3D}
	\end{equation}
	where $S$ is the convex hull of a shape of constant width $K$ and its symmetric with respect to the origin $-K$.
\end{thm} 

It is well known that the inradius of a constant width body $K$ is maximal when it contains a regular tetrahedron \cite[Section 14.3]{bodies_of_constant_width}. Thus $S \subset B(\sqrt{3/8})$. Moreover, like in the previous section, $S = \conv\{K,-K\}$ contains the ball $B(1/2)$, assuming $K$ and $-K$ have the same circumcenter. Denoting again by $\mathcal A$ the family of convex shapes $K$, symmetric with respect to the origin, such that $B(1/2)\subset K \subset B(\sqrt{3/8})$, the Blaschke-Lebesgue theorem could be relaxed to 
\begin{equation}\label{eq:relax3D}
\max_{S\in \mathcal A} (V(B,B,S) - V(B,S,S)).
\end{equation}
However, compared to the two dimensional case, the study of this problem is not obvious. It can be shown that the inradius of $S$ must be equal to $1/2$. The derivative of the mapping
\[ \lambda \mapsto V(B,B,\lambda S) - V(B,\lambda S,\lambda S)\]
at $\lambda = 1$ equals $V(B,B,S)-2V(B,S,S)$. Since $B \subset 2S$, this derivative is negative. Therefore, contractions of $S$ increase the objective function. 

If the Meissner conjecture is true, which according to \cite{kawohl-webe} and simulations presented in \cite{AntunesBogosel22} is quite probable, solution to problem \eqref{eq:max_equiv_3D} is given by the convex hull of one of the Meissner bodies $M$ and its opposite. Assuming $M$ and $-M$ have the same center (i.e., their tetrahedral skeletons $T$, $-T$ have the same center) it can be observed that $S_M=\conv\{M,-M\}$ verifies the following:
\begin{itemize}[noitemsep]
	\item the boundary of $S_M$ can be partitioned into six parts, supported on a square (face of the cube $\conv \{T,-T\}$). 
	\item the six regions are isometric and have two planes of symmetry given by the diagonal planes of the square. 
	\item changing between the two Meissner tetrahedra rotates all six regions in $S_M$ with an angle of $\pi/2$. Thus, the convex hull of $M$ and $-M$ has the same volume, surface area and mean width for the two Meissner tetrahedra.
\end{itemize}

\begin{figure}
	\centering
	\begin{tabular}{cc}
		\includegraphics[height=0.3\textwidth]{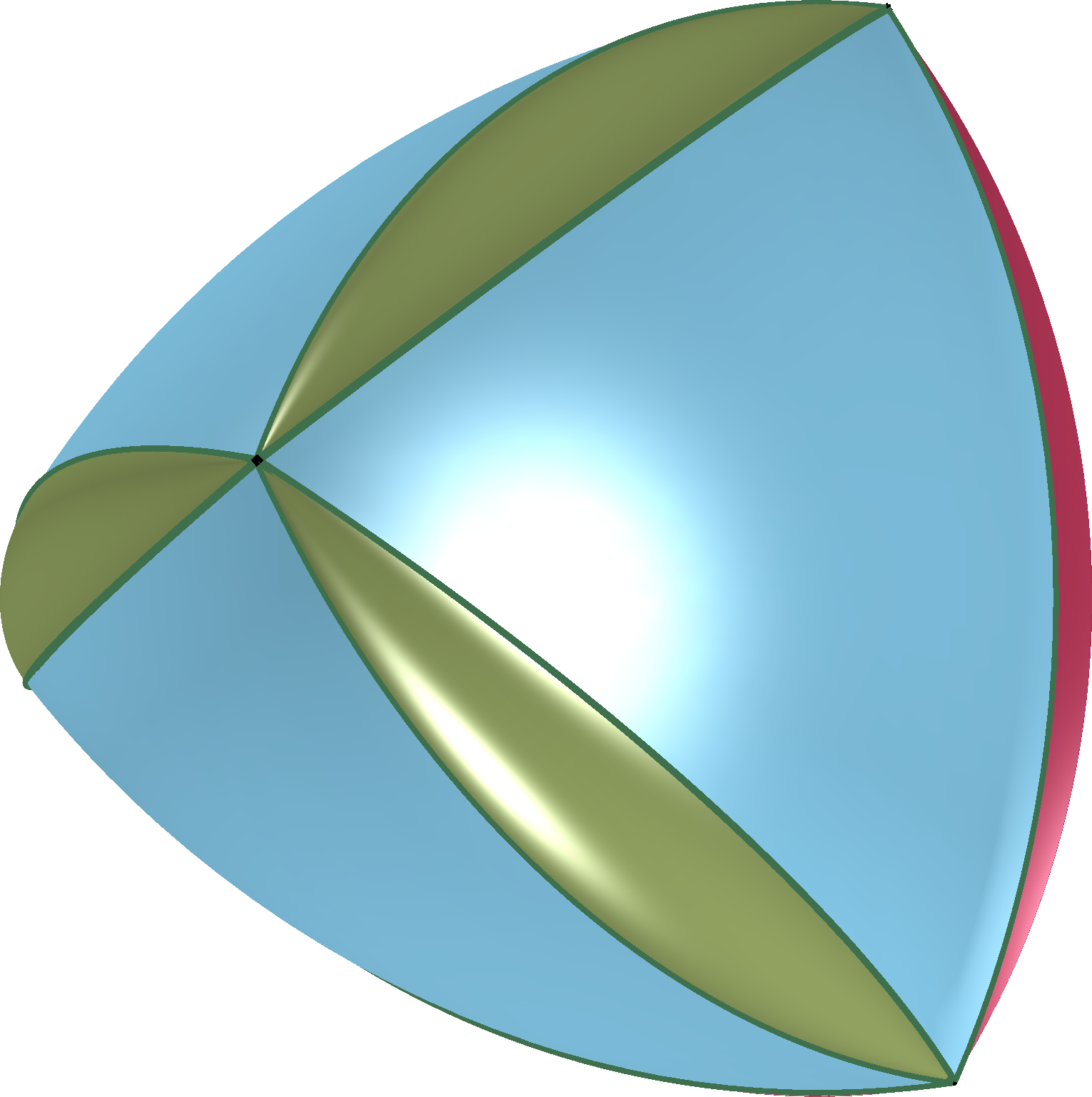} &
		\includegraphics[height=0.3\textwidth]{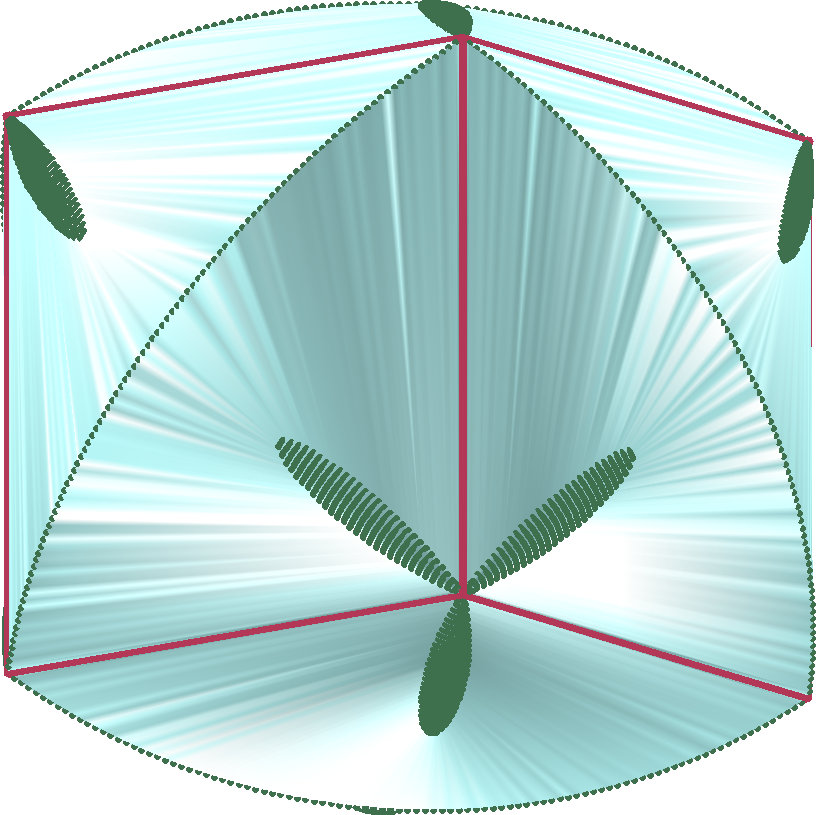} \\
		\includegraphics[height=0.3\textwidth]{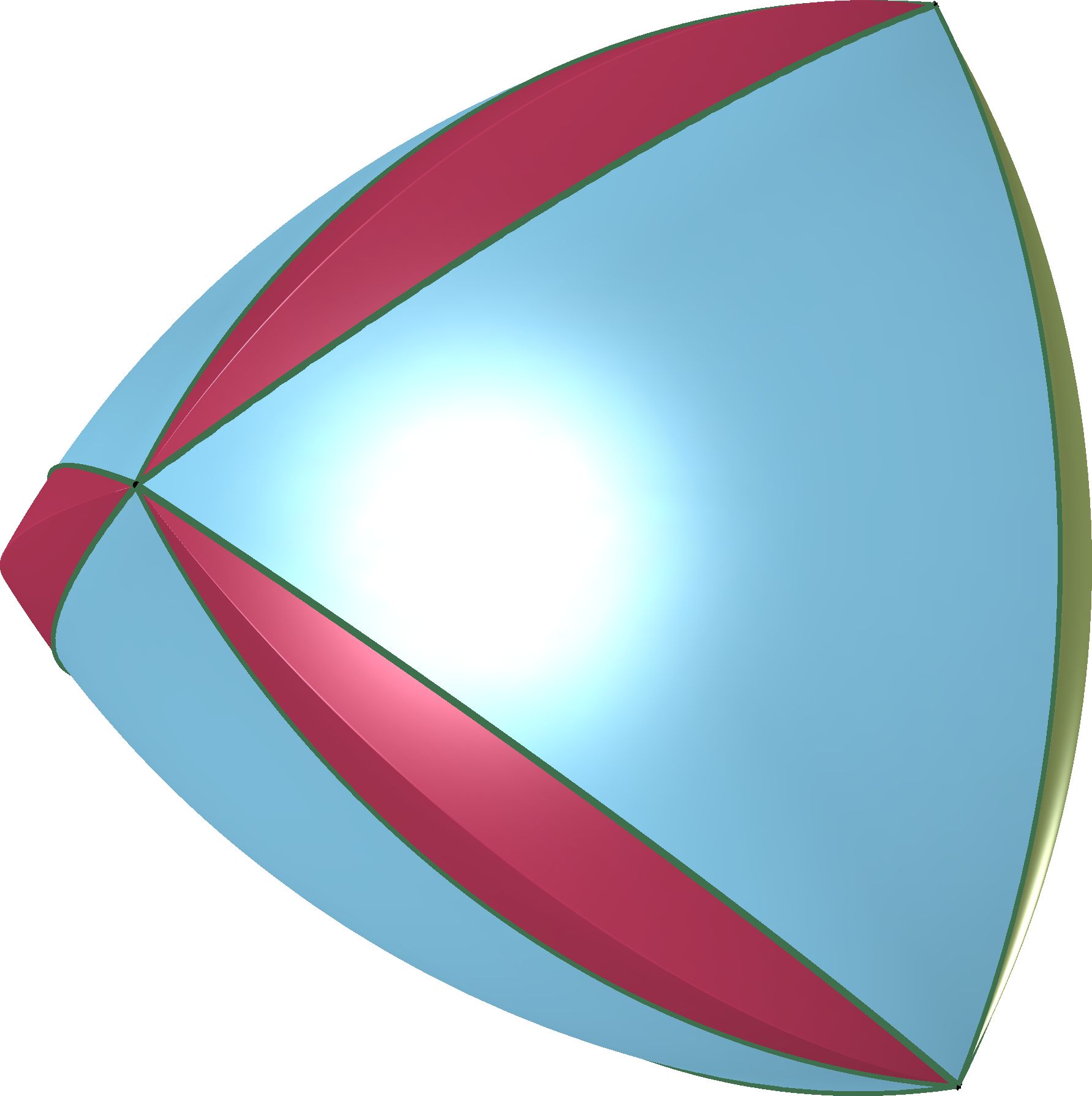} &
		\includegraphics[height=0.3\textwidth]{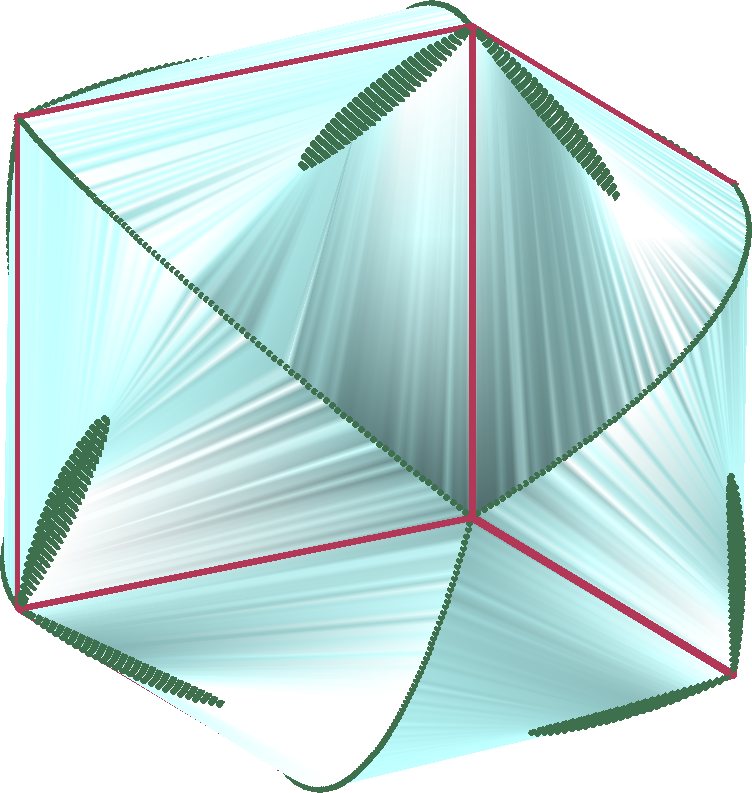} \\
		$M$ & $\conv\{M,-M\}$ 
	\end{tabular}
\caption{Meissner tetrahedra $M$ of the two types together with the corresponding convex hull of $M$ and $-M$, the conjectured solution for \eqref{eq:max_equiv_3D} and \eqref{eq:relax3D}. The faces of $\conv\{M,-M\}$ are changed with a rotation of $\pi/2$ when passing from one tetrahedra to the other one. Red lines on the right column represent the edges of the cube $\conv\{T,-T\}$.}
\label{fig:Meissner}
\end{figure}

In Figure \ref{fig:Meissner} an illustration of $S_M$ is given for the two types of Meissner tetrahedra. Recall that in each of the three pairs of opposite edge, one of them is smoothed. This results in two configurations: either three edges with a common vertex are smoothed or three edges common to the same face are smoothed. In the plot of $S_M$ the red lines represent the edges of the cube $\conv\{T,-T\}$ and the green points are points in $M\cup -M$ which are also on the boundary of the convex hull $\conv\{M,-M\}$. 

The parts of the surface of $S_M$ which are not in common with $M$ or $-M$ are unions of line segments. This is similar to other three dimensional solutions of shape optimization problems under convexity constraint, like those obtained for Newton's problem of least resistance (see \cite{Wachsmuth2014} and the references therein). 

\section{Conclusion}

Mixed volumes are used in dimension two to obtain two proofs of the Blaschke-Lebesgue theorem. Chakerian's proof \cite{Chakerian} was illustrated, showing, in particular, that the mixed volume $A(K,-K)$ when $K$ is a Reuleaux polygon, is in fact constant for any convex set contained between two particular polygons, one of which is the skeleton $\mathcal S(K)$ of $K$ (Theorem \ref{thm:mixed-constant}).

Mixed volumes also allow to reformulate the Blaschke-Lebesgue problem as a reverse isoperimetric problem (Theorem \ref{thm:reformulation}). The equivalent problem admits a relaxation which is solved by the regular hexagon of unit width, according to \cite{bianchini_henrot}, providing a new proof of the Blaschke-Lebesgue problem, which may lead to a extension to higher dimensions.

The proof of Chakerian does not extend to the three dimensional case, since the best known symmetric universal cover, the rhombic dodecahedron, does not give a tight enough estimate for the maximum of the mixed volume $V(K,K,-K)$. On the other hand, the three dimensional Blaschke-Lebesgue problem can also be reformulated as a reverse isoperimetric problem involving the mean width and the surface area (Theorem \ref{thm:reformulation3D}).

While the three dimensional Blaschke-Lebesgue problem remains unsolved, it is possible that the reformulation given in Theorem \ref{thm:reformulation3D} may give a new path to a solution. We conclude with a series of \bo{open questions}:

\begin{enumerate}[noitemsep]
	\item What is the minimal volume (or area) of a centrally symmetric universal cover in dimension three? In dimension two the answer is given by the regular hexagon of unit width, while in dimension three the rhombic dodecahedron is a universal cover. 
	\item Is \eqref{eq:relax3D} a relaxation of \eqref{eq:max_equiv_3D}? In other words, the constraint $S=\conv\{K,-K\}$, where $K,-K$ are constant width bodies with the same circumcenter, can be relaxed to $B(1/2)\subset S \subset B(\sqrt{3/8})$?
	\item Can the results of \cite{bianchini_henrot} be extended to the three dimensional case to characterize solutions of \eqref{eq:relax3D}? 
	\item Problems \eqref{eq:max_equiv_2D_relax} and \eqref{eq:relax3D} have the same structure: one of the terms is linear in terms of the support function, while the other one is quadratic. Moreover, $V(B,B,S)$ is half of the derivative of $t\mapsto V(B,S+tB,S+tB)$ at $t=0$, showing that the two problems have a very particular structure. Can they be solved using a unified framework?
\end{enumerate}

{\bf Acknowledgments.} The author thanks Prof. Rolf Schneider for clarifying aspects related to \eqref{eq:schneider} and \cite{Schneider_mix_vol}. This work was supported by the ANR Shapo program (ANR-18-CE40-0013).

\bibliographystyle{abbrv}
\bibliography{./biblio}

\bigskip
\small\noindent
Beniamin \textsc{Bogosel}: Centre de Math\'ematiques Appliqu\'ees, CNRS,\\
\'Ecole polytechnique, Institut Polytechnique de Paris,\\
91120 Palaiseau, France \\
{\tt beniamin.bogosel@polytechnique.edu}\\
{\tt \nolinkurl{http://www.cmap.polytechnique.fr/~beniamin.bogosel/}}

\end{document}